\documentclass{amsart}
\usepackage{amsmath}
\def\MM{{\mathfrak{M}}}
\def\rr{{\mathfrak{r}}}
\def\JJ{{\mathcal{J}}}
\def\RR{{\mathcal{R}}}
\def\AA{{\mathcal{A}}}\def\GL{\operatorname{GL}}

\newtheorem{theorem}{Theorem}[section]
\newtheorem{lemma}[theorem]{Lemma}

\makeatletter
  
  \@addtoreset{equation}{section}
 \makeatother
\begin{document}
%\date{Version 3a last changed 19 November 2008}
\title[Riemannian Geometric Realizations of algebraic curvature tensors]{RIEMANNIAN GEOMETRIC REALIZATIONS FOR RICCI TENSORS OF
GENERALIZED ALGEBRAIC CURVATURE OPERATORS}

\author{P. Gilkey}
\address{Mathematics Department, University of Oregon\\
Eugene OR 97403 USA\\
E-mail: gilkey@uoregon.edu}

\author{S. Nik\v cevi\'c}

\address{Mathematical Institute, Sanu,
Knez Mihailova 35, p.p. 367\\
11001 Belgrade,
Serbia.\\
E-mail: stanan@mi.sanu.ac.rs}

\author{D. Westerman}
\address{Mathematics Department, University of Oregon\\
Eugene OR 97403 USA\\
E-mail: dwesterm@uoregon.edu}

\begin{abstract}
We examine questions of geometric realizability for algebraic structures which arise naturally in affine and Riemannian
geometry.
\end{abstract}


%\bodymatter
\maketitle
\section{Introduction}\label{sect-1}

Many questions in Riemannian geometry involve constructing
geometric realizations of algebraic objects where the objects in question are invariant under the action of the structure
group $G$. We present several
examples to illustrate this point. We first review previously known results. Section \ref{sect-1.1} deals with Riemannian
algebraic curvature tensors, Section \ref{sect-1.2} deals with Osserman tensors, and Section \ref{sect-1.3} deals
with generalized algebraic curvature operators. 

In Section \ref{sect-1.4} we present the new results of this paper
that deal with a mixture of affine and Riemannian geometry; this mixture has not been considered previously. The results in the
real analytic context can perhaps be considered as extensions of previous results in affine geometry; the results in
the $C^s$ context are genuinely new and require additional estimates. We refer to Section \ref{sect-1.4} for further details.
To simplify the discussion, we shall assume that the underlying dimension $m$ is at least $3$ as the $2$-dimensional case is a
bit exceptional. We adopt the {\it Einstein convention} and sum over repeated indices henceforth.

\subsection{Realizing Riemannian algebraic curvature tensors}\label{sect-1.1} Let $V$ be an $m$-dimensional real vector space
and let
$\rr(V)\subset\otimes^4V^*$ be the set of all {\it Riemannian algebraic
curvature tensors}; $A\in\rr(V)$ if and only if $A$ has the
symmetries of the Riemannian curvature tensor of the Levi-Civita
connection:
\begin{equation}\label{eqn-1}
\begin{array}{l}
A(x,y,z,w)=-A(y,x,z,w),\qquad A(x,y,z,w)=A(z,w,x,y),\\
A(x,y,z,w)+A(y,z,x,w)+A(z,x,y,w)=0\,.\\

\end{array}\end{equation}
Let $A\in\mathfrak{r}(V)$ and let $\langle\cdot,\cdot\rangle$ be a non-degenerate symmetric bilinear form on $V$ of signature
$(p,q)$. The triple $\MM:=(V,\langle\cdot,\cdot\rangle,A)$ is said to be a {\it pseudo-Riemannian algebraic curvature model};
let $\Xi(V)$ be the set of such models.

Let $\mathcal{M}:=(M,g)$ be a pseudo-Riemannian manifold. Let $\nabla^g$ be the associated Levi-Civita connection and let
$R_P^g\in\otimes^4T^*_PM$ be the curvature tensor at a point $P$ of $ M$. Since $R_P^g$ satisfies the symmetries of Equation
(\ref{eqn-1}),  $$\MM_P(\mathcal{M}):=(T_PM,g_P,R_P^g)\in \Xi(T_PM)\,.$$

The following result shows every $\mathfrak{M}\in\Xi(V)$ is geometrically realizable; in particular, the
symmetries of Equation (\ref{eqn-1}) generate the universal symmetries of the curvature tensor of the Levi-Civita connection.
\begin{theorem}\label{thm-1.1}
Let $\MM\in\Xi(V)$. There exists
a pseudo-Riemannian manifold $\mathcal{M}$, a point $P\in M$, and an isomorphism
$\phi$ from $T_PM$ to $V$ so that $\MM_P(\mathcal{M})=\phi^*\MM$.
\end{theorem}

\subsection{Osserman geometry}\label{sect-1.2} 

The
relevant structure group which arises in this context is the orthogonal group
$O(V,\langle\cdot,\cdot\rangle)$; one can ask geometric realization questions concerning any
$O(V,\langle\cdot,\cdot\rangle)$ invariant subset of
$\rr(V)$. If $\MM=(V,\langle\cdot,\cdot\rangle,A)\in\Xi(V)$, the  {\it Jacobi operator}
$\JJ_{\MM}\in\operatorname{End}(V)\otimes V^*$ is characterized by the relation:
$$\langle\mathcal{J}_{\MM}(x)y,z\rangle=A(y,x,x,z)\,.$$

If $p>0$, then $\MM$ is said to be
{\it timelike Osserman} if the spectrum of $\JJ_{\MM}$ is constant on the pseudo-sphere
of unit timelike vectors in $V$. The notion {\it spacelike Osserman} is defined similarly if $q>0$.
If $p>0$ and if $q>0$, work of  N. Bla\v zi\'c et al.\cite{BBG97} and of
Garc\'{\i}a-R\'{\i}o et al.\cite{GKV97} shows these two notions are equivalent and thus we shall simply say $\MM$ is 
{\it Osserman} in this context. As this definition is invariant under the action of the structure group
$O(V,\langle\cdot,\cdot\rangle)$, it extends to the geometric setting. Thus a pseudo-Riemannian manifold $\mathcal{M}$ will be
said to be {\it Osserman} provided that the associated model
$\MM_P(\mathcal{M})$ is Osserman for every $P\in M$. 

Work of Chi \cite{C88} shows there
are
$4$-dimensional Osserman Riemannian algebraic curvature tensors which are not geometrically realizable by Osserman
manifolds. The field is a vast one and we refer to Nikolayevsky\cite{N05} for further details in the Riemannian
setting and to Garc\'{\i}a-R\'{\i}o et al.\cite{GKV02} for a discussion in the pseudo-Riemannian setting; it is
possible to construct many examples of Osserman tensors in the algebraic context which have no corresponding
geometrical analogues.

\subsection{Affine geometry}\label{sect-1.3}

Let $\nabla$ be a torsion free connection on $M$. The associated curvature
operator $\RR\in T^*M\otimes T^*M\otimes\operatorname{End}(TM)$ is a $(3,1)$ tensor which has
the symmetries
\begin{equation}\label{eqn-2}
\RR(x,y)z=-\RR(y,x)z,\quad
\RR(x,y)z+\RR(y,z)x+\RR(z,x)y=0\,.\vphantom{\vrule height 12pt}
\end{equation}
As we are in the affine setting, there is no analogue of the additional curvature symmetry $A(x,y,z,w)=A(z,w,x,y)$
which appears in the pseudo-Riemannian setting. In the algebraic context, let
$\mathfrak{A}(V)\subset V^*\otimes V^*\otimes\operatorname{End}(V)$ be the set of $(3,1)$ tensors satisfying the relations of
Equation (\ref{eqn-2}). An element $\AA\in\mathfrak{A}(V)$ is said to be a {\it generalized algebraic curvature
operator}.

If $\nabla$ is a torsion free connection and if $P\in M$,
then $\RR^\nabla_P\in\mathfrak{A}(T_PM)$. The
following geometric realizability result is closely related to Theorem \ref{thm-1.1}. It shows that any universal
symmetry of the curvature tensor of an affine connection is generated by the summetries of Equation (\ref{eqn-2}).

\begin{theorem}\label{thm-1.2}
Let $\AA\in\mathfrak{A}(V)$. There exists a torsion free connection $\nabla$ on a smooth manifold $M$, a point $P\in M$, and an
isomorphism $\phi$ from $T_PM$ to $V$ so that
$\RR^\nabla_P=\phi^*\AA$.
\end{theorem}

We contract indices to define the {\it Ricci tensor} $\rho(\mathcal{A})\in V^*\otimes V^*$ by setting
$$\rho(\mathcal{A})(x,y):=\operatorname{Trace}\{z\rightarrow\AA(z,x)y\}\,.$$
The decomposition $V^*\otimes V^*=\Lambda^2(V^*)\oplus S^2(V^*)$ sets $\rho(\AA)=\rho_a(\AA)+\rho_s(\AA)$
where $\rho_a(\mathcal{A})$ and $\rho_s(\AA)$ are the antisymmetric and symmetric Ricci tensors. The natural structure group
for $\mathfrak{A}(V)$ is the general linear group $\GL(V)$. The Ricci tensor defines a $\GL(V)$
equivariant short exact sequence
$$
0\rightarrow\ker(\rho)\rightarrow\mathfrak{A}(V)\rightarrow V^*\otimes V^*\rightarrow 0\,.
$$
Strichartz
\cite{S88} showed this short exact sequence is $\GL(V)$ equivariantly split and gives a
$\GL(V)$ equivariant decomposition
$$\mathfrak{A}(V)=\ker(\rho)\oplus\Lambda^2(V^*)\oplus S^2(V^*)$$
into irreducible $\GL(V)$ modules. The {\it Weyl projective
curvature operator} $\mathcal{P}(\AA)$ is the projection of $\AA$ on $\ker(\rho)$;
$\AA$ is said to be {\it projectively flat} if $\mathcal{P}(\AA)=0$, 
$\AA$ is said to be {\it Ricci symmetric} if $\rho_a(\AA)=0$, and $\AA$ is said to be {\it Ricci
antisymmetric} if $\rho_s(\AA)=0$. These notions for a connection are defined similarly. There are 8 natural geometric
realization questions which arise in this context and whose realizability\cite{GNW08} may be summarized in the following table
-- the possibly non-zero components being indicated by $\star$:
$$\begin{array}{|c|c|c|r||c|c|c|r|}\noalign{\hrule}
\ker(\rho)&S^2(V^*)&\Lambda^2(V^*)&&\ker(\rho)&S^2(V^*)&\Lambda^2(V^*)&\\
\noalign{\hrule}\star&\star&\star&\text{yes}&0&\star&\star&\text{yes}\\
\noalign{\hrule}\star&\star&0&\text{yes}&0&\star&0&\text{yes}\\
\noalign{\hrule}\star&0&\star&\text{yes}&0&0&\star&\text{no}\\
\noalign{\hrule}\star&0&0&\text{yes}&0&0&0&\text{yes}\\
\noalign{\hrule}
\end{array}$$
Thus, for example, if $\AA$ is projectively flat and Ricci
symmetric, then $\AA$ can be geometrically realized by a projectively flat Ricci symmetric
torsion free connection. But if $\AA\ne0$ is projectively flat and Ricci antisymmetric,
then $\AA$ can not be geometrically realized by a projectively flat Ricci antisymmetric torsion free connection.

\subsection{Torsion free connections and Riemannian geometry}\label{sect-1.4}

We now combine the settings of Sections \ref{sect-1.1} and \ref{sect-1.3}. Let
$\langle\cdot,\cdot\rangle$ be a non-degenerate symmetric inner product on $V$ of signature $(p,q)$. Fix a basis $\{e_i\}$ for $V$ and let $g_{ij}:=\langle e_i,e_j\rangle$ give the
components of $\langle\cdot,\cdot\rangle$. Let $g^{ij}$ be the inverse matrix. If $\mathcal{A}\in\mathfrak{A}(V)$, expand
$\AA(e_i,e_j)e_k=\AA_{ijk}{}^\ell e_\ell$. The scalar curvature
$\tau$ and trace free Ricci tensor are then given, respectively, by
$$
\tau(\AA,\langle\cdot,\cdot\rangle):=g^{ij}\AA_{kij}{}^k,\quad
\rho_0(\AA,\langle\cdot,\cdot\rangle):=
\rho_s(\AA)-{\textstyle\frac{\tau(\mathcal{A},\langle\cdot,\cdot\rangle)}m}\langle\cdot,\cdot\rangle\,.
$$
Let $S_0^2(V^*,\langle\cdot,\cdot\rangle)$ be the space of trace free symmetric bilinear forms. One has an
$O(V,\langle\cdot,\cdot\rangle)$ invariant decomposition of $V^*\otimes V^*$ into irreducible
$O(V,\langle\cdot,\cdot\rangle)$ modules
$$
V^*\otimes V^*=\Lambda^2(V^*)\oplus S_0^2(V^*,\langle\cdot,\cdot\rangle)\oplus\mathbb{R}\,.
$$
This decomposition leads to 8 geometric realization questions which are natural with respect
to the structure group
$O(V,\langle\cdot,\cdot\rangle)$ and which can all be solved either in the real analytic category or in the $C^s$ category of
$s$-times differentiability for any $s\ge1$. The following is the main result of this paper; as our considerations are
local, we take $M=V$ and $P=0$.

\begin{theorem}\label{thm-1.3}
Let $g$ be a $C^s$ (resp. real analytic) pseudo-Riemannian metric on $V$. Let $\AA\in\mathfrak{A}(V)$. There exists a torsion
free
$C^s$ (resp. real analytic) connection $\nabla$ defined on a neighborhood of $0$ in $V$such that:
\begin{enumerate}
\item $\mathcal{R}^\nabla_0=\AA$.
\item $\nabla$ has constant scalar curvature.
\item If $\AA$ is Ricci symmetric, then $\nabla$ is Ricci symmetric.
\item If $\AA$ is Ricci antisymmetric, then $\nabla$ is Ricci antisymmetric.
\item If $\AA$ is Ricci traceless, then $\nabla$ is Ricci traceless.
\end{enumerate}\end{theorem}

The subspace $\ker(\rho)\subset\mathfrak{A}(V)$ is not an irreducible $O(V,\langle\cdot,\cdot\rangle)$ module but
decomposes as the direct sum of 5 additional irreducible factors -- see Bokan\cite{B90}. This decomposition will play no role
in our further discussion and studying the additional realization questions which arise from this decomposition
is a topic for future investigation.

\section{The proof of Theorem \ref{thm-1.3}}\label{sect-2}

We assume $s\ge1$ and $m\ge3$ henceforth; fix
$\mathcal{A}\in\mathfrak{A}(V)$. Introduce the following notational conventions. Choose a basis
$\{e_i\}$ for
$V$ to identify
$M=V=\mathbb{R}^m$ and let
$\{x_1,...,x_m\}$ be the associated coordinates.

For $\delta>0$, let $B_\delta:=\{x\in\mathbb{R}^m:|x|<\delta\}$ where $|x|$ is the
usual Euclidean norm on $\mathbb{R}^m$. Let $C^s_\delta$ be the set of functions on $B_\delta$ which are
$s$-times differentiable. Let
$\alpha=(\alpha_1,...,\alpha_m)$ be a multi-index. Set
$$
\partial_i:={\textstyle\frac{\partial}{\partial x_i}},\quad
\partial_x^\alpha:=(\partial_1)^{\alpha_1}...(\partial_m)^{\alpha_m},
\quad|\alpha|=\alpha_1+...+\alpha_m\,.$$
If $\mathfrak{Z}$ is a real vector space, let $C^s_\delta(\mathfrak{Z})$ be the set of $C^s$ functions on $B_\delta$ with values
in $\mathfrak{Z}$.  Fix a basis $\{f_\sigma\}$
for $\mathfrak{Z}$ and expand
$P\in C^s_\delta(\mathfrak{Z})$ as
$P=P^\sigma f_\sigma$ for $P^\sigma\in C^s_\delta$. Let $\nu\in\mathbb{R}$. Set
$$|P|:=\sup_\sigma|P^\sigma|\in C^0_\delta\quad\text{and}\quad
||P||_{\delta,\nu,-1}:=0\,.$$
For $0\le r\le s$, define
$||P||_{\delta,\nu,r}\in[0,\infty]$ by setting
$$||P||_{\delta,\nu,r}:=\sup_{|\alpha|=r,\ |x|<\delta}|\partial_x^\alpha P(x)|\cdot|x|^{-\nu}\,.$$
Thus $||P||_{\delta,\nu,r}\le C$ implies $|\partial_x^\alpha P(x)|\le C|x|^\nu$ for $|\alpha|=r$ and $|x|<\delta$.  Let
$$\mathfrak{G}:=S^2((\mathbb{R}^m)^*)\otimes\mathbb{R}^m,\quad
\mathfrak{S}:=S^2((\mathbb{R}^m)^*),\quad\text{and}\quad\mathfrak{A}:=\mathfrak{A}(\mathbb{R}^m)\,.$$
We use the basis $\{e_i\}$ and the coordinate
frame
$\{\partial_i\}$ to determine the components of tensors of all types; if computing relative to some orthonormal frame $\{E_i\}$,
we shall make this explicit. Thus, for example, if
$\mathcal{S}\in\mathfrak{S}$, then $\mathcal{S}_{ij}=\mathcal{S}(e_i,e_j)$ while if $\mathcal{S}\in C^s_\delta(\mathfrak{S})$,
then
$\mathcal{S}_{ij}:=\mathcal{S}(\partial_{x_i},\partial_{x_j})$.
If
$\Gamma,\mathcal{E}\in C^s_\delta(\mathfrak{G})$, define $\mathcal{L}(\Gamma)\in C^{s-1}(\mathfrak{A})$ and
$\Gamma\star\mathcal{E}\in C^s_\delta(\mathfrak{A})$ by setting
\begin{equation}\label{eqn-3}
\begin{array}{l}
\mathcal{L}(\Gamma)_{ijk}{}^l:=\partial_i\Gamma_{jk}{}^l-\partial_j\Gamma_{ik}{}^l,\\
(\Gamma\star\mathcal{E})_{ijk}{}^\ell:=\mathcal{E}_{in}{}^\ell\Gamma_{jk}{}^n+
\Gamma_{in}{}^\ell\mathcal{E}_{jk}{}^n-\mathcal{E}_{jn}{}^\ell\Gamma_{ik}{}^n-\Gamma_{jn}{}^\ell\mathcal{E}_{ik}{}^n
\,.\phantom{\vrule height 12pt}
\end{array}\end{equation}
 If $\Gamma\in C^s_\delta(\mathfrak{G})$, let $\nabla(\Gamma)$ be the $C^s$
torsion free connection on $B_\delta$ with Christoffel symbol $\Gamma$. One has:
\begin{equation}\label{eqn-4}\begin{array}{l}
\mathcal{R}^{\nabla(\Gamma)}=\mathcal{L}(\Gamma)+\textstyle\frac12\Gamma\star\Gamma,\\
\rho(\Gamma\star\Gamma)_{jk}=2\Gamma_{\ell n}{}^\ell\Gamma_{jk}{}^n-2\Gamma_{jn}{}^\ell\Gamma_{\ell k}{}^n
=\rho(\Gamma\star\Gamma)_{kj},\phantom{\vrule height 10pt}\\
\rho_a(\mathcal{R}^{\nabla(\Gamma)})_{jk}=\rho_a(\mathcal{L}(\Gamma))_{jk}=\textstyle\frac12\left\{\partial_k\Gamma_{ji}{}^i
-\partial_j\Gamma_{ki}{}^i\right\}\,.\phantom{\vrule height 10pt}
\end{array}\end{equation}
One says that $\Gamma$ is {\it
normalized} if
\medbreak\qquad(1) $\Gamma(0)=0$ and $\RR^{\nabla(\Gamma)}=\AA+O(|x|^2)$.
\par\qquad(2) $\rho_s(\mathcal{R}^{\nabla(\Gamma)})$ is $C^s$.
\par\qquad(3) $\rho_a(\mathcal{R}^{\nabla(\Gamma)})(\partial_i,\partial_j)=\rho_a(\mathcal{A})(e_i,e_j)$ on $B_\delta$.
\medbreak\noindent
We remark that Assertion (2) is non-trivial as $\mathcal{R}^{\nabla(\Gamma)}$ need only be $C^{s-1}$. This is a technical
condition used subsequently to avoid loss of smoothness.

 Theorem \ref{thm-1.2} follows from the following observation which forms
the starting point in our proof of Theorem
\ref{thm-1.3}:
\begin{lemma}\label{lem-2.1}
If $\Gamma_{uv}{}^l:=\textstyle\frac13(\AA_{wuv}{}^l+\AA_{wvu}{}^l)x^w$, then $\Gamma$ is normalized.
\end{lemma}
\begin{proof}
Since $\Gamma(0)=0$, one has:
\begin{eqnarray*}
&&\mathcal{R}^{\nabla(\Gamma)}_0(\partial_i,\partial_j)\partial_k
=\textstyle\left\{\partial_i\Gamma_{jk}{}^l(0)-\partial_j\Gamma_{ik}{}^l(0)\right\}\partial_\ell\\
&&=\textstyle\frac13\left\{\AA_{ijk}{}^l+\AA_{ikj}{}^l-\AA_{jik}{}^l-\AA_{jki}{}^l\right\}\partial_\ell\\
&&=\textstyle\frac13\left\{\AA_{ijk}{}^l-\AA_{kij}{}^l+\AA_{ijk}{}^l-\AA_{jki}{}^l\right\}\partial_\ell
=\textstyle \AA_{ijk}{}^l\partial_\ell\,.
\end{eqnarray*}
By Equation (\ref{eqn-4}), $\rho_a(\mathcal{R}^{\nabla(\Gamma)})_{ij}=\rho_a(\mathcal{L}(\Gamma))_{ij}
=\rho_a(\mathcal{R}^{\nabla(\Gamma)}_0)_{ij}=\rho_a(\mathcal{A})_{ij}$.
\end{proof}

We continue our analysis with the following basic solvability result:

\begin{lemma}\label{lem-2.2} If $\Theta\in C^s_\delta(\mathfrak{S})$, then there exists $\mathcal{E}\in
C^s_\delta(\mathfrak{G})$ so
$\rho(\mathcal{L}(\mathcal{E}))=\Theta$, so $\mathcal{E}_{ij}{}^j=0$, and so
$||\mathcal{E}||_{\delta,\nu+1,r}\le||\Theta||_{\delta,\nu,r}+r||\Theta||_{\delta,\nu+1,r-1}$.
\end{lemma}

\begin{proof} By assumption, $m\ge3$. For each pair of indices $\{i,j\}$, not necessarily distinct,
choose $k=k(i,j)=k(j,i)$ with $k\ne i$ and $k\ne j$. Set
$$
\mathcal{E}_{ij}{}^\ell:=\left\{\begin{array}{lll}
\textstyle \int_0^{x_k}\Theta_{ij}(x_1,...,x_{k-1},u,x_{k+1},...,x_m)du&\text{ if }&\ell=k,\\
0&\text{ if }&\ell\ne k\,.
\end{array}\right.$$
Since $k\ne j$, $\mathcal{E}_{ij}{}^j=0$. Consequently Equation (\ref{eqn-3}) yields
$$\rho(\mathcal{L}(\mathcal{E}))_{ij}=\partial_\ell\mathcal{E}_{ij}{}^\ell=\Theta_{ij}\,.$$
Expand $\partial_x^\alpha=\partial_k^\mu\partial_x^\beta$ where $\beta$ does not involve the index $k$. Then:
$$
\partial_x^\alpha\mathcal{E}_{ij}{}^k=\left\{\begin{array}{lll}
\textstyle\int_0^{x_{k}}\partial_x^\alpha\Theta_{ij}(x_1,...,x_{k-1},u,x_{k+1},...,x_m)du&\text{if}&\mu=0,\\
\mu\partial_k^{\mu-1}\partial_x^\beta\Theta_{ij}&\text{if}&\mu>0\,.
\end{array}\right.$$
Assume that $|\partial_x^\alpha\Theta_{ij}(x)|\le C|x|^\nu$ for all $x\in B_\delta$ and all $|\alpha|=j$. Then
\begin{eqnarray*}
&&\textstyle|\int_0^{x_k}\partial_x^\alpha\Theta_{ij}(x_1,...,x_{k-1},u,x_{k+1},...,x_m)du|\\
&\le&\textstyle|x_k|\int_0^1|\partial_x^\alpha\Theta_{ij}(x_1,...,x_{k-1},tx_k,x_{k+1},...,x_m)|dt\\
&\le&|x_k|\cdot C|x|^\nu\le C|x|^{\nu+1}\,.
\end{eqnarray*}
The estimates of the Lemma now follow.
\end{proof}

Let $g$ be a $C^s$ pseudo-Riemannian metric on $B_\delta$ for $\delta<1$, let $\{E_i\}$ be a $C^s$ $g$-orthonormal frame for the
tangent bundle of $B_\delta$, and let
$e_i:=E_i(0)$. Let $\Gamma\in C^s_\delta(\mathfrak{G})$.  Define $\Theta=\Theta(\Gamma)\in
C^s_\delta(\mathfrak{S})$ by:
\begin{equation}\label{eqn-5}
\Theta_{ij}:=\rho_s(\mathcal{R}^{\nabla(\Gamma)})(E_i,E_j)-\rho_s(\AA)(e_i,e_j)\,.
\end{equation}
Use Lemma \ref{lem-2.2} to define
$\mathcal{E}=\mathcal{E}(\Gamma)\in C^s_\delta(\mathfrak{G})$ so that
$\rho_s(\mathcal{L}(\mathcal{E}))=-\Theta$.
We use Lemma \ref{lem-2.1} to choose an initial Christoffel symbol $\Gamma_1\in C^s_\delta(\mathfrak{G})$ which is normalized.
Inductively, set 
$$\Theta_\nu:=\Theta(\Gamma_\nu),\quad
\mathcal{E}_{\nu+1}:=\mathcal{E}(\Gamma_\nu),\quad
\Gamma_{\nu+1}:=\Gamma_\nu+\mathcal{E}_{\nu+1}\,.$$
We will set $\Gamma_\infty:=\Gamma_1+\mathcal{E}_2+...$, we will establish convergence, and we will show $\Gamma_\infty$
defines a connection with the desired properties. We begin by using Equation (\ref{eqn-4}) to compute:
\begin{eqnarray*}
&&\Theta_{\nu+1,ij}=\rho_s(\mathcal{R}_\nu)(E_i,E_j)-\rho_s(\AA)(e_i,e_j)+\rho_s(\mathcal{L}(\mathcal{E}_{\nu+1}))_{ij}\\
&+&\rho_s(\mathcal{L}(\mathcal{E}_{\nu+1}))(E_i,E_j)-\rho_s(\mathcal{L}(\mathcal{E}_{\nu+1}))_{ij}\\
&+&\rho_s\{(\Gamma_\nu+\textstyle\frac12\mathcal{E}_{\nu+1})\star\mathcal{E}_{\nu+1}\}(E_i,E_j)\,.
\end{eqnarray*}
As $\rho_s(\mathcal{L}(\mathcal{E}_{\nu+1}))=-\Theta_\nu$, the first line vanishes and
\begin{equation}\label{eqn-6}
\Theta_{\nu+1,ij}=-\Theta_\nu(E_i,E_j)+\Theta_{\nu,ij}
+\rho_s\{(\Gamma_\nu+\textstyle\frac12\mathcal{E}_{\nu+1})\star\mathcal{E}_{\nu+1}\}(E_i,E_j)\,.
\end{equation}
Choose $\kappa\ge1$ so we have the following estimates for any $x\in B_\delta$:
\begin{eqnarray*}
&&|\mathcal{S}(E_i,E_j)-\mathcal{S}_{ij}|\le\kappa|\mathcal{S}|\cdot|x|^2\quad\forall\quad\mathcal{S}\in
C^0_\delta(\mathfrak{S}),\\
&&|\rho_s(\Gamma\star\mathcal{E})(E_i,E_j)|\le\kappa\{|\Gamma|\cdot|\mathcal{E}|\}\quad\forall\quad\Gamma,\mathcal{E}\in
C^0_\delta(\mathfrak{G})\,.
\end{eqnarray*}

\begin{lemma}\label{lem-2.3}
Adopt the notation established above. Then $\Gamma_\nu$ is normalized for all $\nu$. Furthermore, there exists $\delta_0>0$ and
there exist constants
$C_r>0$  for $0\le r\le s$ so for $\nu=1,2,...$ we have the estimates:
\begin{enumerate}
\item $||\Gamma_\nu||_{\delta_0,1-r,r}\le\frac1{4\kappa}C_r$.
\item $||\Theta_\nu||_{\delta_0,2\nu-r,r}\le C_r^\nu$.
\item $||\mathcal{E}_{\nu+1}||_{\delta_0,2\nu+1-r,r}\le C_r^\nu+rC_{r-1}^\nu$.
\end{enumerate}
\end{lemma}

\begin{proof} By assumption $\Gamma_1$ is normalized. We assume inductively $\Gamma_\nu$ is normalized and show
 $\Gamma_{\nu+1}$ is normalized. As
$\mathcal{E}_{\nu+1,ij}{}^j=0$, Equation (\ref{eqn-4}) yields
\begin{equation}\label{eqn-7}
\rho_a(\mathcal{R}_{\nu+1})_{ij}=\rho_a(\mathcal{R}_\nu)_{ij}=\rho_a(\mathcal{A})_{ij}\quad
\text{on}\quad B_\delta\,.
\end{equation}
Since $\Theta_\nu=\rho_s(\mathcal{R}_\nu)(E_i,E_j)-\rho_s(\mathcal{A})_{ij}=O(|x|^2)$, $\mathcal{E}_{\nu+1}=O(|x|^3)$ and
$$\mathcal{R}_{\nu+1}=\mathcal{R}_\nu+O(|x|^2)=\mathcal{A}+O(|x|^2)\,.$$
As $\Gamma_\nu$ is normalized, $\Theta_\nu\in C_\delta^s(\mathfrak{S})$. Hence $\mathcal{E}_{\nu+1}\in C_\delta^s(\mathfrak{G})$
and
$\Gamma_{\nu+1}\in C_\delta^s(\mathfrak{G})$. Since $\rho_s(\mathcal{L}(\mathcal{E}_{\nu+1}))=-\Theta_\nu$ is $C^s$,
we may conclude that
$\rho_s(\mathcal{R}_{\nu+1})$ is $C^s$ even though $\mathcal{R}_{\nu+1}$ need only be $C^{s-1}$. Thus
$\Gamma_{\nu+1}$ is normalized.

We establish the estimates by induction on $r$ and then on $\nu$; Assertion (3)${}_{\nu,r}$ follows from Assertions
(2)${}_{\nu,r}$ and (2)${}_{\nu,r-1}$ and from Lemma
\ref{lem-2.2}. Suppose first that $r=0$; this is, somewhat surprisingly, the most difficult case. As $\Gamma_1$ is
normalized,
$\mathcal{R}_1=\AA+O(|x|^2)$. One has
$E_i(0)=\partial_i$. Thus
$\Theta_1=O(|x|^2)$. As $\Gamma_1=O(|x|)$, by shrinking $\delta$, we may choose $\bar C_0$ so
$$|\Gamma_1|(x)\le \bar C_0|x|\quad\text{and}\quad|\Theta_1|(x)\le \bar C_0|x|^2\quad\text{on}\quad B_\delta\,.$$
Choose $C_0$ and $\delta_0<\delta<1$ so that
$$
\bar C_0+1<\textstyle\frac1{4\kappa}C_0<C_0,
\quad\kappa+\textstyle\frac14C_0+\textstyle\frac12\kappa\le C_0,\quad
\delta_0^2C_0<\frac12\,.
$$
If $\nu=1$, then Assertions (1) and (2) follow from the choices made. Assume the Assertions hold for
$\mu\le\nu$ where $\nu\ge1$. Then
\begin{eqnarray*}
|\Gamma_{\nu+1}|&\le&|\Gamma_1|+|\mathcal{E}_2|+...+|\mathcal{E}_{\nu+1}|
\le\bar C_0|x|+C_0|x|^3+C_0^2|x|^5+...\\
&\le&\bar C_0|x|+|x|\textstyle\frac{C_0|x|^2}{1-C_0|x|^2}\le (\bar C_0+1)|x|\le\frac1{4\kappa}C_0\,.
\end{eqnarray*}
We use Equation (\ref{eqn-6}) to complete the induction step for $r=0$ by checking
\begin{eqnarray*}
|\Theta_{\nu+1}|&\le&\kappa\{|x|^2\cdot|\Theta_\nu|+(|\Gamma_\nu|+|\mathcal{E}_{\nu+1}|)|\mathcal{E}_{\nu+1}|\}\\
&\le&\kappa\{C_0^\nu|x|^{2\nu+2}+(\textstyle\frac1{4\kappa}C_0|x|+C_0^\nu|x|^{2\nu+1})C_0^\nu|x|^{2\nu+1}\}\\
&\le&C_0^\nu|x|^{2\nu+2}\{\kappa+\textstyle\frac14C_0+\textstyle\frac12\kappa\}\le C_0^{\nu+1}|x|^{2\nu+2}\,.
\end{eqnarray*}

We now suppose $r=1$; we get 1 less power of $|x|$ in
the decay estimates. We choose $\bar C_1$ so $|\partial_k\Gamma_1|\le\bar C_1$ and $|\partial_k\Theta_1|\le\bar C_1|x|$; the
desired estimates then hold for $\nu=1$ for $C_1$ sufficiently large. We proceed by induction on $\nu$. We then have for
sufficiently large $C_1$ and small $\delta_0$ that:
\begin{eqnarray*}
&&|\partial_k\Gamma_{\nu+1}|\le|\partial_k\Gamma_1|+|\partial_k\mathcal{E}_2|+...+|\partial_k\mathcal{E}_{\nu+1}|\\
&\le&\bar C_1+\{C_1+C_0\}|x|^2+\{C_1^2+C_0^2\}|x|^4+...\le\textstyle\frac1{4\kappa}C_1\,.
\end{eqnarray*}
We differentiate Equation (\ref{eqn-6}) to obtain
\begin{eqnarray*}
&&\partial_k\Theta_{\nu+1,ij}=-(\partial_k\Theta_\nu)(E_i,E_j)+\partial_k\Theta_{\nu,ij}-\Theta_\nu(\partial_kE_i,E_j)
-\Theta_\nu(E_i,\partial_kE_j)\\
&&\quad+\rho_s\{(\partial_k\Gamma_\nu+\textstyle\frac12\partial_k\mathcal{E}_{\nu+1})\star\mathcal{E}_{\nu+1}
+(\Gamma_\nu+\textstyle\frac12\mathcal{E}_{\nu+1})\star\partial_k\mathcal{E}_{\nu+1}\}(E_i,E_j)\\
&&\quad+\rho_s\{(\Gamma_\nu+\textstyle\frac12\mathcal{E}_{\nu+1})\star\mathcal{E}_{\nu+1}\}\{(\partial_kE_i,E_j)
+(E_i,\partial_kE_j)\}\,.
\end{eqnarray*}
Thus for a suitably chosen constant $\kappa_1=\kappa_1(\mathcal{A},E,\Gamma_1)$ which is independent of $\nu$ and for suitably
chosen
$C_1>C_0$, we have
\begin{eqnarray*}
&&|\partial_k\Theta_{\nu+1}|\le\kappa_1\{|x|^2|\partial_k\Theta_\nu|+|x|\cdot|\Theta_\nu|
+(|\partial_k\Gamma_{\nu}|+|\partial_k\mathcal{E}_{\nu+1}|)|\mathcal{E}_{\nu+1}|\\
&&\quad+(|\Gamma_{\nu}|+|\mathcal{E}_{\nu+1}|)\cdot|\partial_k\mathcal{E}_{\nu+1}|
+(|\Gamma_\nu|+|\mathcal{E}_{\nu+1}|)\cdot|\mathcal{E}_{\nu+1}|\}\\
&&\quad\le\textstyle\kappa_1|x|^{2\nu+1}\{C_1^\nu+C_0^\nu+(\frac1{4\kappa}C_1+(C_1^\nu+C_0^\nu)|x|^{2\nu})C_0^\nu\\
&&\quad+(\textstyle\frac1{4\kappa}C_0+C_0^\nu|x|^{2\nu})(C_1^\nu+C_0^\nu)+(\frac1{4\kappa}C_0|x|+C_0^\nu|x|^{2\nu+1})C_0^\nu\}\,.
\end{eqnarray*}
A crucial point is that
there are no $C_1^{\nu+1}$ terms present. The desired estimate now follows for $C_1$ sufficiently large and $\delta_0$
sufficiently small. This completes the proof of the case
$r=1$; the higher order derivatives are estimated similarly.
\end{proof}

Since $|\partial_x^\alpha\mathcal{E}_{\nu+1}|\le C_r^{\nu}|x|^{2\nu+1-r}$, the series $\mathcal{E}_2+\mathcal{E}_3+...$
converges geometrically for small $x$ and thus the sequence $\Gamma_\nu$ converges in the
$C^r$ topology to a limit $\Gamma_\infty$. Note that as we have to shrink $\delta$ at each stage, we do not get convergence in
the
$C^\infty$ topology even if the initial metric is smooth. We use Equation (\ref{eqn-5}) to see that for
small $x$ we have:
\begin{equation}\label{eqn-8}
\rho_a(\mathcal{R}_\infty(x))_{ij}=\lim_{\nu\rightarrow\infty}
\rho_a(\mathcal{R}_\nu(x))_{ij}=\rho_a(\AA)_{ij}\,.
\end{equation}
This controls the antisymmetric part of the Ricci tensor. To control the symmetric part of the Ricci tensor, we use the
$g$-orthonormal frame $\{E_i\}$. We compute, using Equation (\ref{eqn-7}), that:
\begin{equation}\label{eqn-9}
\begin{array}{l}
\rho_s(\mathcal{R}_\infty(x))(E_i,E_j)=\lim_{\nu\rightarrow\infty}\rho_s(\mathcal{R}_\nu(x))(E_i,E_j)\\
\quad=\lim_{\nu\rightarrow\infty}\Theta_{\nu,ij}(x)+\rho_s(\AA)_{ij}=\rho_s(\AA)_{ij}\,.\vphantom{\vrule height 11pt}
\end{array}\end{equation}
The frame $\{E_i\}$ is $g$-orthonormal. Thus $\mathcal{R}_\infty$ has constant scalar curvature. By Equation (\ref{eqn-8}) if
$\AA$ is Ricci symmetric, then so is $\mathcal{R}_\infty$. By Equation (\ref{eqn-9}), if $\AA$ is Ricci antisymmetric or is
Ricci tracefree, so is
$\mathcal{R}_\infty$. This completes the proof of Theorem \ref{thm-1.3} in the $C^s$ category.

In the real analytic category, we complexify and consider the complex ball of radius $\delta$ in $\mathbb{C}^m$. Since
 $C^0$ convergence of holomorphic functions gives convergence in the holomorphic setting, Theorem \ref{thm-1.3}
follows in the real analytic context as well.

\section*{Acknowledgments} Research of P. Gilkey supported by Project MTM2006-01432 (Spain) and PIP 6303-2006-2008 Conicet
(Argentina). Research of S. Nik\v cevi\'c supported by Project 144032 (Srbija). Research of D. Westerman supported by the
University of Oregon.

\end{document}